\documentclass[superscriptaddress,twocolumn,notitlepage,pra,longbibliography]{revtex4-2}
\usepackage[colorlinks=true,citecolor=blue,urlcolor=black]{hyperref}
\usepackage{graphicx}
\usepackage{dcolumn}
\usepackage{bm}
\usepackage{amsthm}
\usepackage{amsmath}
\usepackage{amssymb}
\usepackage{color}
\usepackage{multirow}
\usepackage{algorithm}
\usepackage{algorithmicx}
\usepackage{algpseudocode}
\usepackage{subfigure}
\usepackage{hyperref}
\usepackage{wasysym}
\usepackage{xcolor}
\usepackage{booktabs}

\usepackage{float}
\usepackage[toc,page]{appendix}
\usepackage{comment}

\usepackage{threeparttable}

\usepackage{longtable}

\definecolor{Red}{rgb}{1,0,0}

\begin{document}

\title{Quantum-Inspired Mean Field Probabilistic Model for Combinatorial Optimization Problems}

%

\author{Yuhan Huang}
\affiliation{Department of Electronic and Computer Engineering, The Hong Kong University of Science and Technology, 999077, Hong Kong}

\author{Siyuan Jin}
\affiliation{Department of information systems, The Hong Kong University of Science and Technology, 999077, Hong Kong}

\author{Yichi Zhang}
\affiliation{Department of Physics, The Hong Kong University of Science and Technology, 999077, Hong Kong}

\author{Ling Pan}
\affiliation{Department of Electronic and Computer Engineering, The Hong Kong University of Science and Technology, 999077, Hong Kong}

\author{Qiming Shao}
\email{eeqshao@ust.hk}
\affiliation{Department of Electronic and Computer Engineering, The Hong Kong University of Science and Technology, 999077, Hong Kong}
\affiliation{IAS Center for Quantum Technologies, The Hong Kong University of Science and Technology, 999077, Hong Kong}


\begin{abstract}
Combinatorial optimization problems are pivotal across many fields. Among these, Quadratic Unconstrained Binary Optimization (QUBO) problems, central to fields like portfolio optimization, network design, and computational biology, are NP-hard and require exponential computational resources. To address these challenges, we develop a novel Quantum-Inspired Mean Field (QIMF) probabilistic model that approximates solutions to QUBO problems with enhanced accuracy and efficiency. The QIMF model draws inspiration from quantum measurement principles and leverages the mean field probabilistic model. We incorporate a measurement grouping technique and an amplitude-based shot allocation strategy, both critical for optimizing cost functions with a polynomial speedup over traditional methods. Our extensive empirical studies demonstrate significant improvements in solution evaluation for large-scale problems of portfolio selection, the weighted maxcut problem, and the Ising model. Specifically, using S\&P 500 data from 2022 and 2023, QIMF improves cost values by 152.8\% and 12.5\%, respectively, compared to the state-of-the-art baselines. Furthermore, when evaluated on increasingly larger datasets for QUBO problems, QIMF's scalability demonstrates its potential for large-scale QUBO challenges.
\end{abstract}

\maketitle
\section{Introduction}

Quadratic Unconstrained Binary Optimization (QUBO) problems have gained significant attention due to their ability to model a wide range of complex scenarios across various domains, including finance~\cite{lang2022strategic}, operation research~\cite{barahona1988application}, and physics~\cite{lucas2014ising}. The objective of QUBO problems is to find a set of binary variables that minimizes or maximizes a quadratic polynomial (typically represented as a sum of linear and quadratic terms)~\cite{kochenberger2014unconstrained}. This work defines QUBO problem as $x'Vx$, where \(x \in \{0, 1\}^n\) and \(V\) is a \(n \times n\) symmetric matrix. This problem is a fundamental problem in graph theory, where the binary variables $x$ represent the selection or rejection of the vertices, and the matrix $V$ encodes the interaction weights between pairs of vertices. 

Traditional exact solving methods, such as cutting plane algorithms~\cite{de1994cutting} and branch-and-bound techniques~\cite{rendl2007branch}, can indeed solve QUBO problems to optimality. However, as problem size grows, these methods face exponential computational complexity, limiting their use for large-scale QUBO problems.
To mitigate this issue, various approximate solving methodologies have been developed, including simulated annealing~\cite{kagawa2020fully}, genetic algorithms~\cite{kim2019comparison}, and probabilistic modeling techniques~\cite{batsheva2024protes}. Finding nearly optimal solutions involves iterative searches to assess solution quality using a cost function, where lower scores indicate better solutions. However, the increasing computational complexity of the cost function for larger problem scales remains a significant challenge, even for approximate methods. As problem size increases, evaluating the cost function becomes too time-consuming, reducing the efficiency of approximate methods.


Quantum computing harnesses the unique properties of quantum mechanics, such as superposition and entanglement, allowing multiple potential solutions simultaneously, and offering potential advantages over classical algorithms~\cite{nielsen2001quantum} (a more detailed description is outlined in Appendix.~\ref{QC}). The QUBO problem can be mapped into the Ising model, where the objective function is the energy function. In this setup, finding the optimal solution to a QUBO problem equates to determining the state with the minimum energy. A detailed discussion of the quantum algorithms for QUBO is provided in Appendix.~\ref{QA_QUBO}. Although quantum annealing~\cite{pastorello2019quantum} and the Quantum Approximate Optimization Algorithm~\cite{farhi2014quantum} demonstrate the potential for significantly enhanced optimization capabilities two significant challenges persist: i) current hardware does not satisfy the requirements of quantum algorithms, and ii) the quantum measurement process, which is essential for extracting information from a quantum system, remains a time-consuming bottleneck in quantum computing.

\begin{figure*}[hptb]
  \centering
  \includegraphics[width=0.8\textwidth]{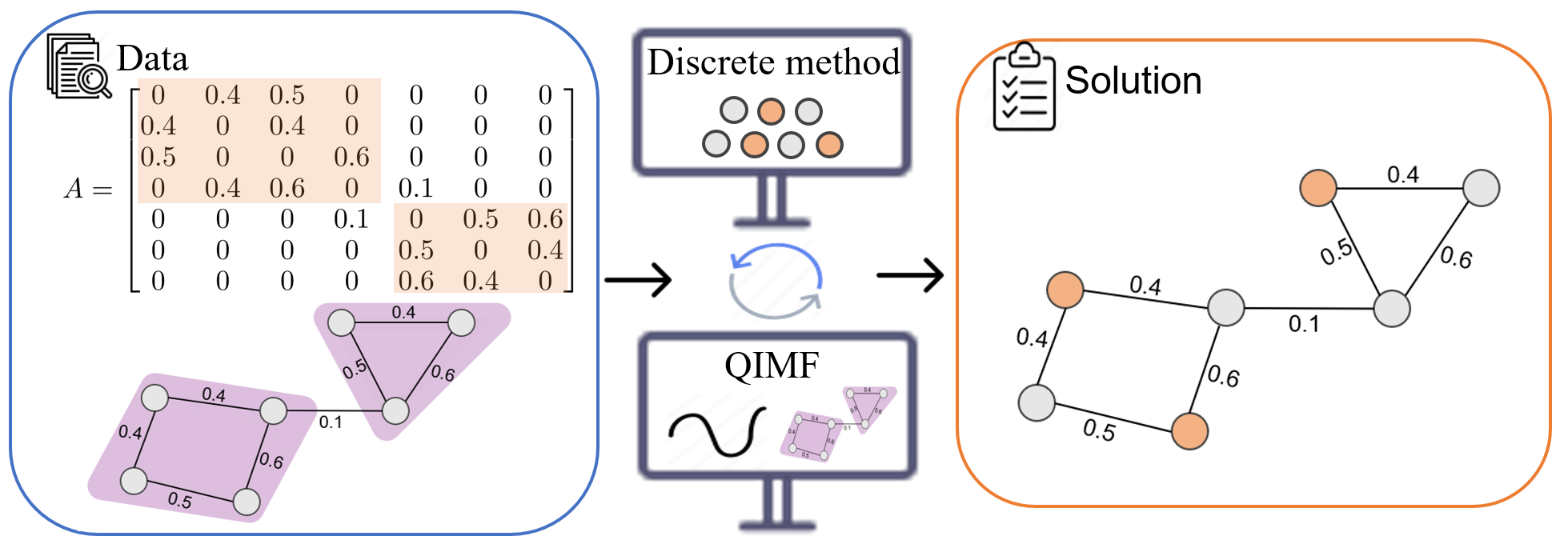}  
  \caption{The QIMF method for solving QUBO problems, utilizing a continuous model and incorporating block information, with depicted output results on the right.}
  \label{fig:qubo} 
\end{figure*}

This paper proposes a quantum-inspired method for solving QUBO problems, specifically the Quantum inspired mean field (QIMF) probabilistic model, as illustrated in Fig.~\ref{fig:qubo}. We adopt a mean field probabilistic model~\cite{zhang2022differentiable} while leveraging a quantum-inspired technique~\cite{verteletskyi2020measurement, yen2020measuring,zhu2023optimizing} to efficiently compute the cost function. QIMF is designed to have the following
features: \textbf{i) continuous optimization transformation:} Mean field probabilistic model within the QIMF framework converts discrete optimization problems (binary selections of 0 or 1) into continuous domains using a softmax function. This method utilizes the simplest probabilistic models, referred to as independent category probabilistic models in energy modeling. By converting the problems into continuous domains, the method enables the use of gradient-based optimizers, which efficiently find optimal solutions through derivatives. \textbf{ii) acceleration of our method:} Using quantum-inspired methods, we implement measurement grouping~\cite{verteletskyi2020measurement, yen2020measuring} and amplitude-based shot allocation~\cite{zhu2023optimizing} to enhance computational efficiency. Measurement grouping enables parallel computation of each term in a mapped Ising model, while amplitude shot allocation adjusts computational resources based on the coefficients of each term. These techniques collectively expedite the cost function calculation for datasets with block-like structures and non-uniform distributions. These characteristics reflect the topology of the underlying graph, manifesting as clusters with varying edge weights.
\textbf{iii) real-world data suitability:} QIMF is particularly effective with data that naturally fits the method's partitioning requirements, such as the S\&P 500 stock data, which is characterized by block-like structures and non-uniform distribution.

The main contributions of our research are outlined as follows:

\begin{itemize}
    \item We have developed the QIMF method, which facilitates the continuous solution of discrete optimization problems and enhances the speed of the evaluation process.
    \item Our approach is theoretically proven to accelerate solutions for QUBO problems characterized by non-uniform and block-like structured data.
    \item Experimental results validate the superiority of our method, particularly when applied to real-world financial data, demonstrating significant improvements over existing techniques.
\end{itemize}



\section{Related Works}
\label{related-works}

\textbf{WSBM: Weighted Stochastic Block Model.}
The WSBM extends the Stochastic Block Models (SBM) by incorporating weights into the traditional framework, making it a more general form of SBM tailored for analyzing both the presence and strength of interactions within networks~\cite{airoldi2008mixed}. The model is defined as WSBM($\mathbf{n}$, $N$, $P$, $W$), where a graph's nodes are organized into \(N\) blocks, with \(\mathbf{n} = (n_1, n_2, \dots, n_N)\) representing the vector of node counts in each block. The interactions between these blocks are characterized by a probabilistic matrix \(P\) and further detailed through a weight matrix \(W\), providing a comprehensive framework for analyzing structured network data.

WSBM is successfully applied to many real-world problems, including social network analysis~\cite{decelle2011asymptotic,airoldi2008mixed, karrer2011stochastic}, protein function prediction~\cite{stanley2019stochastic}, and communication networks~\cite{decelle2011asymptotic}, as well as human brain studies~\cite{faskowitz2018weighted} and in recommendation systems~\cite{xiao2019hybrid}. Setting necessary conditions is the key to optimal adaptive sampling, achieving near-perfect clustering in the WSBM framework~\cite{yun2019optimal}. Additionally, entropy calculations play a critical role in WSBM research, helping to better identify and understand block structures and elucidate the roles of topological features in empirical networks~\cite{peixoto2012entropy}.

\textbf{QUBO: Quadratic Unconstrained Binary Optimization Problem.}
QUBO problem involves finding a binary vector that minimizes a given quadratic polynomial, a task known to be NP-hard and emblematic of combinatorial optimization difficulties. Traditional computational methods include utilizing analytical methods~\cite{caprara1999exact, de1994cutting} and heuristic approaches.  Notably, evolutionary algorithms~\cite{samorani2019clustering} and simulated annealing~\cite{kagawa2020fully} offer approximate yet frequently satisfactory solutions by exploring the search space iteratively to find near-optimal solutions. 
In QUBO problems, data can be represented using the WSBM format (detailed in Appendix~\ref{WSBM}), characterized by diagonal blocks, connectivity matrices, and weight distribution matrices, which we call WSBM-QUBO data. For instance, consider solving $x^{\dagger}Vx$, where $V$ is defined by the vector $\mathbf{n}$, the number of blocks $N$, the connectivity matrix $P$, and the weight distribution matrix $W$.

\textbf{Probabilistic Model.}
Probabilistic models~\cite{zhang2022differentiable} are mathematical tools that express the uncertainty within systems by assigning probabilities to different possible outcomes. These models are particularly powerful in addressing discrete optimization problems, where decisions must be made under uncertainty, and the outcomes are discrete in nature. By capturing the probability distribution over possible solutions, probabilistic models facilitate the exploration of the solution space more effectively, allowing for the identification of near-optimal solutions following the optimization of the model. In our research, we concentrate on utilizing the most fundamental form of these models, known as the mean field model, where each variable is treated as independent.

\textbf{Quantum Measurement Optimization.}
Measurement processes are optimized to efficiently extract information from quantum systems, particularly when dealing with the Hamiltonian or energy operator \( H \). The concept of qubit-wise commutativity (QWC)~\cite{verteletskyi2020measurement, yen2020measuring,izmaylov2019unitary}, which is detailed in Appendix~\ref{QWC}, plays a central role in this optimization. By grouping terms of the Hamiltonian into QWC groups, each of which can be measured simultaneously, the process significantly reduces the complexity and the number of separate measurement setups required. 
Another pivotal method is the shot allocation strategy~\cite{zhu2023optimizing}, which further reduces measurement complexity by reallocating the number of shots among measurement operators in the variational quantum algorithm. In variational quantum eigensolver, the energy expectation value  $\langle H \rangle$ is given by the sum of individual terms weighted by their coefficients $\langle\sum_{m} a_m H_m \rangle$. The amplitude-based shot strategy allocates shots to measurement operator terms based on their amplitudes $a_m$. In this strategy, terms with larger amplitudes receive more shots, as they contribute more significantly to the energy estimation. By assigning shots in proportion to the amplitudes, the amplitude-based shot strategy leads to faster convergence and requires fewer overall shots, making it an attractive choice in practical implementations.





\section{QIMF: Quantum Inspired Mean Field Probabilistic Model}
\label{sec:main-method}

This section introduces the QIMF Method, which includes two primary components: the quantum-inspired cost function, and the mean field probabilistic model.

\subsection{Quantum Inspired Cost Function}
\label{subsec:quantum-cost-function}

This subsection discusses the cost functions used in the QUBO problem. We begin by discussing both the classical and quantum formulations of the cost function in QUBO scenarios. Subsequently, we introduce a dequantized version of the cost function. 




\textbf{Classical Cost Function in QUBO.}
The classical cost function, indicated as $\text{Cost}(x)$, is a quadratic form represented by the equation
\begin{align}
    \text{Cost}(x) = x^{\dagger}Vx, ~~~~~V \in \text{WSBM($\mathbf{n}$,N,P,W)},
\end{align}
where vector x indicates binary decision variables, and matrix V is constructed by the WSBM Model.

\textbf{Quantum Cost Function in QUBO.}
The quantum embedding, in the QUBO problem, typically involves the Ising model, which is expressed as follows
\begin{align}
  H = \sum_{ij}V_{ij}Z_iZ_j +\sum_i V_{ii} Z_i=\sum_m a_m H_m, ~~~~~V \in \text{WSBM($\mathbf{n}$,N,P,W)}.    \label{quantum_cf}
\end{align}
Each term $H_m$ can take the form of either $Z_i$ or $Z_j Z_i$. The commutative property, $[H_{mi}, H_{mj}] = 0$, holds for these terms. The cost function is defined as $\text{Cost}(\theta)=\langle 0|U^{\dagger}(\theta) H U(\theta)|0\rangle$, and $U(\theta)$ represent the parameterized quantum circuit.

\textbf{Quantum Inspired Cost Function in QUBO.}
Before discussing the dequantized cost function, it's important to note that data preprocessing in quantum-inspired computing is optional and its effectiveness varies with data complexity, as detailed in Appendix.~\ref{Prep}. The dequantized cost function is given through the grouping technique, detailed shown in the Appendix.~\ref{QI_cost} $\text{Cost}(x) = \sum_{m=1}^{n_w} a_m \lambda_{mi} I(\phi_i == x)$, where indicator function $I$ evaluates to 1 if the condition $\phi_i == x$ is satisfied, and 0 otherwise, the coefficients $a_m$ represent the magnitudes of the corresponding elements in the matrix $V$. The terms $\lambda_{mi}$ and $\phi_i$ are parameters encoding the variable connectivity, and they correspond to the eigenvalue and eigenstate of $H_m$ in Eq.~\ref{quantum_cf}, respectively.

In scenarios where data are organized into a large number of blocks (\( N \)), to reduce the computational complexity of the cost function, we employ an amplitude-based shot allocation method. This approach involves approximating the cost function by distributing various shots across different terms in the summation. We introduce this shot-based approximate cost function as $\text{Cost}_s(x)$, detailed by
\begin{align}
    \text{Cost}_s(x) = \sum_{s=1}^{n_s} (a_s \lambda_{si} I(\phi_i == x)),
\end{align}
where $n_s$ is significantly lower than $n_w$, which denotes a more efficient computation by considering only a subset of all terms. The selection of the term $a_s \lambda_{si} I(\phi_i==x)$ for each shot is governed by a probability $p \sim \frac{|a_m|^2}{\sum |a_m|^2}$. This probability serves to prioritize the terms with larger magnitudes, thereby approximating the full cost function in a resource-efficient manner.

\subsection{Mean Field Probabilistic Model}

The mean field probabilistic model, denoted as $P(\mathbf{X}, \mathbf{\alpha})$, represents the joint probability distribution of a set of independent and identically distributed (IID) variables $\mathbf{X}$. Applying the mean field approach, the joint probability distribution $P$ can be factorized into the product of individual probabilities $p(\mathbf{X}_i, \mathbf{\alpha}_i)$, with each $\mathbf{X}_i$ following the IID property. The behavior of the individual variable is characterized by the probability distribution $p$, which can be modeled using the softmax function $p(\mathbf{X}_i = j, \alpha_{i}) = \frac{e^{\alpha_{ij}}}{\sum_k e^{\alpha_{ik}}}$. To generate samples from the distribution, Monte Carlo sampling techniques are employed. These samples $x \sim P$ are then evaluated $\text{Cost}(x)$, and an objective function sums up these evaluations
\begin{align}
    \ell = \sum_{x \sim P} \text{Cost}(x) = \sum_x \frac{P(\mathbf{X}, \alpha)}{Z(\alpha)} \text{Cost}(x), \label{Eq: ell}
\end{align}
where normalization factor $Z(\alpha)$ ensures that the probabilities are properly scaled. To optimize the probabilistic model, the Monte Carlo gradient is calculated as follows
\begin{align}
    \nabla_{\alpha} \ell = \sum_{x \sim P} \nabla_{\alpha} \ln P(\mathbf{X}, \alpha) \text{Cost}(x),\label{Eq: delta_ell}
\end{align}
where $\nabla_{\alpha_{ij}} ln P(\mathbf{X}_i=\xi)=I(j==\xi)-P(\mathbf{X}_i=\xi)$. See Appendix~\ref{g_PM} for derivation details.


\subsection{Quantum Inspired Mean Field Probabilistic Model}\label{QIMF}
\begin{figure*}[ht]
  \centering
  \includegraphics[width=1\linewidth]{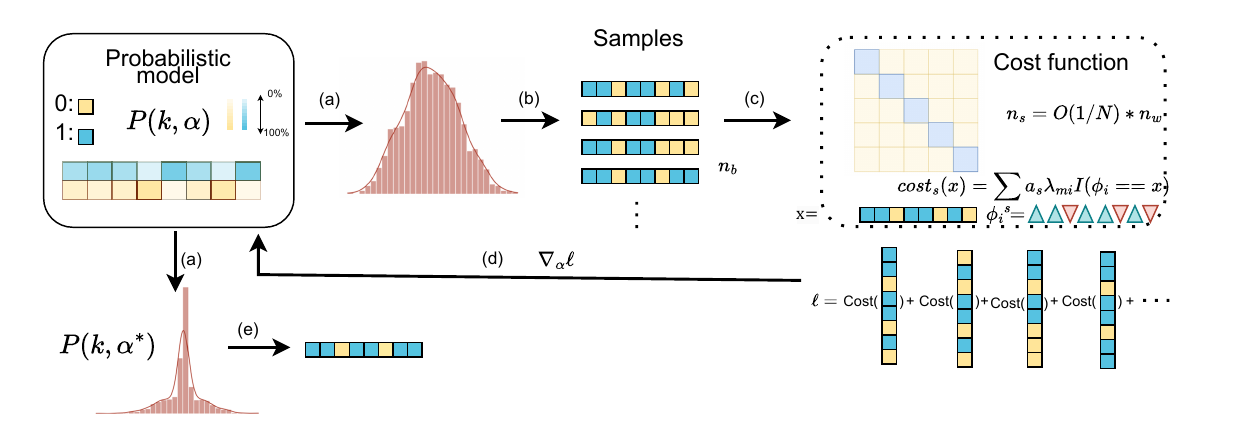}
  \caption{QIMF: Quantum-Inspired Mean Field Probabilistic Model. (a) Representation of the solution distribution through a probabilistic model: Blue for value 1 and yellow for 0, with deeper colors indicating higher probabilities. (b) Sampling of solution candidates from the solution distribution using Monte Carlo methods. (c) Evaluation of the QIMF using the quantum-inspired cost function, where the loss function is assessed via Monte Carlo samples from (b). (d) Optimization of the probabilistic model parameters using a gradient-based method. (e) After convergence, final sampling from the probabilistic model.}
  \label{Fig: QIMF_draw}
\end{figure*}

The QIMF algorithm, shown in Fig.~\ref{Fig: QIMF_draw}, is parameterized by a triplet $(n_b, n_s, n_e)$: the number of candidate solutions $n_b$, the number of shots $n_s$ for evaluating the cost function, and the number of epochs $n_e$ defining the epoch count.

\textbf{1. Probabilistic Initialization.} A probabilistic model $P(\mathbf{X}, \mathbf{\alpha^0)}$, is initialized to represent the search space, where  $\mathbf{X}$ represents a set of IID variables and $\mathbf{\alpha^0}$ includes the initial parameters. Given that QUBO problems are binary, the initial parameter matrix $\mathbf{\alpha^0}$ is structured as an $(\sum \mathbf{n})$ by 2 matrix, where $\sum \mathbf{n}$ represents the number of binary variables to be optimized.

\textbf{2. Monte Carlo Sampling.} A Monte Carlo approach is employed to sample a collection of potential solutions $\{x^i\}_{i=1}^{n_b}$ from the evolving probabilistic model $P(\mathbf{X}, \mathbf{\alpha})$ at each epoch, facilitating the exploration of the solution space.

\textbf{3. Objective Evaluation.} The algorithm assesses the sampled solutions using the objective function $\ell=\sum_{x^i \in P(X, \alpha)}\text{Cost}_s(x^i)$. Following this evaluation, the parameter set $\mathbf{\alpha}$ is updated using the ADAM optimization algorithm, a gradient-based method.

\textbf{4. Gradient Estimation.} The algorithm calculates the gradients of the objective function $\nabla_\alpha \ell=\sum_{x^i \in P(\mathbf{X} \alpha)} \nabla_\alpha \ln P(\mathbf{X}, \mathbf{\alpha}) \text{Cost}_s(x^i)$, leveraging Monte Carlo techniques to estimate the partial derivatives with respect to the parameters $\mathbf{\alpha}$. The gradient of the probabilistic model is determined by $\nabla_{\alpha_{ij}} ln P(X_i=\xi)=I(j==\xi)-P(X_i=\xi)$. 

\textbf{5. Solution Sampling.} Upon converging to the optimal parameters $\mathbf{\alpha^*}$, the probabilistic model $P(\mathbf{X}, \mathbf{\alpha^*})$ becomes capable of generating candidate solutions with the highest quality. Then we can sample the optimal solutions for the QUBO problem. 

\section{Experiments}\label{exper}
We present four cases illustrating the effectiveness of our method. In Case 1, we show empirical evidence of QIMF's polynomial speedup over traditional algorithms. Case 2 compares QIMF with existing optimizers, demonstrating its superiority in various QUBO problem sizes under identical query complexities. Case 3, an ablation study, examines performance changes relative to the number of shots, $n_s$, and the number of samples, $n_b$. Case 4 highlights QIMF's scalability. All experiments run on an AMD EPYC 7542 32-core processor.

\subsection{Polynomial Speedup of QIMF}

Portfolios are divided into \( N \) classes, each containing \( n_i \) portfolios, with internal and external correlations modeled by different normal distributions. This setup illustrates stronger intra-class and weaker inter-class correlations. The portfolio optimization problem is formulated as
\begin{align}
    \arg\min_x \left( \sum_i \lambda V_{ii} x_i^2 + \sum_{i,j} \lambda V_{ij} x_i x_j - \sum_i (1-\lambda) r_i x_i \right),
\end{align} 
where \( V \) is derived from the WSBM(\(\mathbf{n}\), \( N \), \( P \), \( W \)) model with \( N = 5 \) and \( n = [10, 10, 10, 10, 10] \). Connectivity probabilities \( P \) have diagonal elements at 0.2 and off-diagonal at 0.05, and weight distributions \( W \) are \( \mathcal{N}(0, 0.2) \) on the diagonal and \( \mathcal{N}(0, 0.05) \) off-diagonal. Portfolio returns \( r \) follow \( \mathcal{N}(0, 0.2) \), with a balanced risk-return parameter \( \lambda = 0.5 \). The model evaluates problems 1, 2, 3, and 4 under 50 portfolios using two approaches: QUAMF employs a classical cost function, while QIMF utilizes a quantum-inspired cost function. The number of Monte Carlo samples is set to \( n_b = 40 \) with 3k iterations, and the number of shots \( n_s \) is determined by \( \frac{n_w}{N} \).

\begin{figure*}[htbp]
    \centering
    \includegraphics[width=1\textwidth]{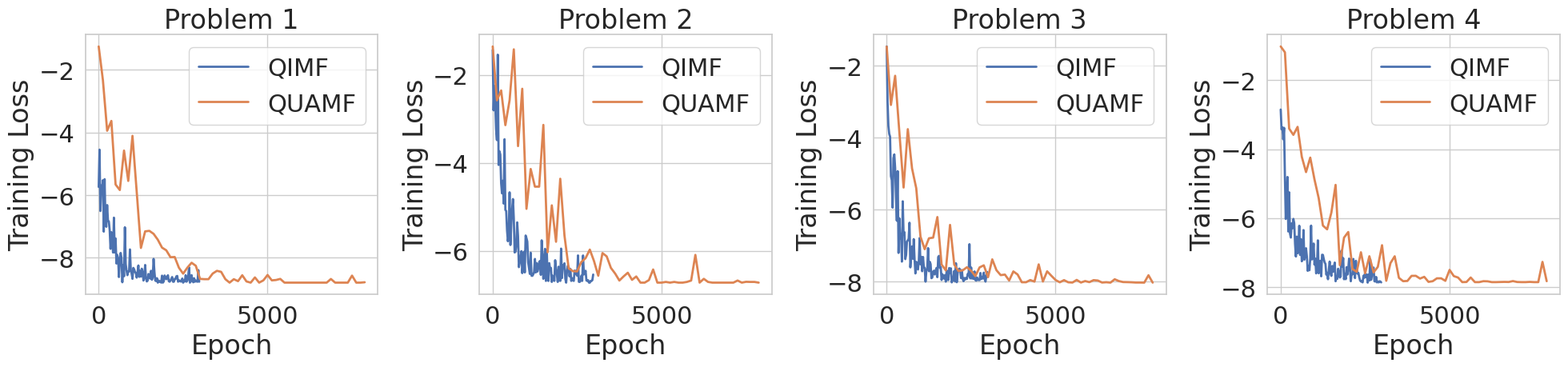} 
    \caption{Simulation Results for four portfolio optimization problems with 50 portfolios; showcasing the comparison of the QIMF, which utilizes a mean field probabilistic model with a quantum-inspired cost function, against the QUAMF, employing a classical cost function for solving QUBO problems.}
    \label{fig:case1}
\end{figure*}
 
Although both QIMF and QUAMF are iterated 3k times, the computational complexity of QUAMF's cost function is N times greater than that of QIMF, which leads to differing scales in iteration in the figures. To accurately reflect the computational effort disparity, we recalibrate using epochs as a unit of measure. Specifically, for QIMF, each epoch corresponds directly to an iteration, whereas for QUAMF, N epochs correspond to one iteration. In the presented results, we truncate the epochs at 8000. The simulation results in Fig.~\ref{fig:case1}, illustrate that our method converges faster than traditional methods by a polynomial factor in different portfolio optimization problems.

\subsection{QIMF versus Existing Optimizers}


In case 2, we employ the QIMF method for addressing portfolio optimization problems~(PO) across different scales, specifically with 50~(PO1), and 500~(PO2, PO3) variables. Additionally, we apply QIMF to the weighted max cut problem~(WM) for graphs with 30~(WM1), 50~(WM2), and 500~(WM3) nodes, as well as to the Ising model~(IS) sized at 50~(IS1) and 500~(IS2). All these applications adhere to the WSBM model with N blocks. Specifically, the classification of companies into 114 distinct industries within the S\&P 500 index data directly influences our choice of the parameter $n_s$. We implement the QIMF method, setting $n_s$ as $\frac{n_w}{N}$, where $N=114$ is the number of industries, for our analyses of the year 2022 in PO2, and for 2023 in PO3~\footnote{S\&P 500 stock data for the years 2022 and 2023 is sourced from~\url{https://www.kaggle.com/datasets/pavankrishnanarne/s-and-p-500-stock-data-from-listing-day-to-2023/data}}. For a more detailed discussion on the determination of the value of N, refer to Appendix.~\ref{PO}.

Our QIMF approach is rigorously benchmarked against six established strategies. These include the commercial optimization solver Gurobi~(BS1)~\cite{gurobi}, Reinforcement Learning~(BS2)~\cite{zhang2023let},  Greedy~(BS3)~\cite{vince2002framework}, OnePlusOne~(BS4)~\cite{bennet2021nevergrad}, Simulated Annealing~(BS5)~\cite{bennet2021nevergrad}, and PROTES~(BS6)~\cite{batsheva2024protes}. Among these, BS2 utilizes unsupervised learning, BS3 through BS5 are heuristic methods, and BS6 incorporates a probabilistic model enhanced by Monte Carlo sampling.

\begin{table*}[htbp]
\centering
\caption{Performance comparison of optimization algorithms}
\label{tab:comparison_results_horizontal}
\scalebox{0.8}{
\begin{tabular}{@{}c c cccc cc cc@{}}
\toprule
\multirow{2}{*}{Method}&\multirow{2}{*}{Type} & \multicolumn{3}{c}{PO} & \multicolumn{3}{c}{WM} & \multicolumn{2}{c}{IS} \\
& & PO1(1e-4) $\downarrow$ & PO2(1e-3) $\downarrow$ & PO3(1e-3) $\downarrow$ & WM1 $\uparrow$ & WM2 $\uparrow$& WM3 $\uparrow$  & IS1 $\downarrow$& IS2 $\downarrow$ \\
\hline
\addlinespace[2pt] 
BS1 & C & $ \textbf{-3.51}$&-&-&$\textbf{42.70}$&\textbf{66.31}&-&\textbf{-12.95}&-\\
\cmidrule(lr){3-5} \cmidrule(lr){6-8} \cmidrule(lr){9-10}
BS2 & UL & -&-&-&$37.51 \pm 0.07$& $57.18 \pm 0.46$ &205.31&-&- \\
\cmidrule(lr){3-5} \cmidrule(lr){6-8} \cmidrule(lr){9-10}
BS3 & H & $ \textbf{-3.51} \pm \textbf{0}$&-2.35 &-5.12&$40.07\pm 0.25$&$60.94\pm 2.55$&265.41&$-12.41\pm0.12$&-13.14\\
BS4 & H &$-3.34 \pm 3.59$e-6&14.49&-1.95&$39.63\pm16.54$&$62.98\pm 8.02$&202.02&$-10.63\pm3.35$&-4.58\\
BS5 & H &$-3.43\pm 6.19$e-7&0.71&-4.32&$39.56\pm4.97$&$63.58\pm 1.98$&204.81&$-11.91\pm0.53$&-7.24\\
\cmidrule(lr){3-5} \cmidrule(lr){6-8} \cmidrule(lr){9-10}
BS6 & MC &$-3.17\pm1.37$e-6&-0.15& -0.85&$25.13\pm0.15$&$30.93\pm0.39$&163.41&$-8.33\pm 0.21$&-7.02 \\
Ours & MC & $ \textbf{-3.51} \pm \textbf{0}$ &\textbf{-5.94}&\textbf{-5.76}&$\text{41.15}\pm\text{1.53}$& $\text{63.88} \pm \text{0.62}$&\textbf{267.47}& $\text{-12.45} \pm \text{0.00} $&\textbf{-15.54}\\
\bottomrule
\end{tabular}}
\footnotesize
\\ 
\begin{tablenotes}
\item
\textit{Note:} "PO" is the portfolio optimization problem, "WM" is the weighted maxcut problem, and "IS" is the Ising problem. For PO, we utilize the real S\&P dataset; the cost value represents the risk minus return. Risk is calculated based on the correlation of return ratios for each portfolio, and return is calculated based on either the annual or monthly returns. For WM, we employ a designed maxcut problem; the cost value represents the maximum cut. For IS, we use a designed Ising problem; the cost value represents the ground energy. The values are presented as mean $\pm$ variance. The arrows $\downarrow$ and $\uparrow$ indicate that lower and higher values, respectively, are preferable. Bold values denote the best performance in each category. 
\end{tablenotes}
\end{table*}

We report the comparative performance of these algorithms, evaluated under identical computational query complexities in Table.~\ref{tab:comparison_results_horizontal}. In scenarios involving small-sized problems, Gurobi achieves the optimal cost score across portfolio optimization, weighted maxcut, and the Ising model. Nevertheless, Gurobi's applicability is constrained by its scalability challenges in handling large-scale applications. Our method, in contrast, consistently exhibits lower risk in portfolio optimization scenarios, successfully achieves the maximum cut in weighted maxcut cases, and identifies the minimal ground energy in the Ising model, outperforming baseline methods BS2, BS3, BS4, BS5, and BS6 and closely matching the performance of BS1. Importantly, our method substantially outperforms all other baseline strategies in large-scale systems such as the PO2, PO3, WM3, and IS2 problems. Specifically, in the PO1, PO2, and PO3 scenarios involving real datasets from the financial market. These results suggest that our approach is exceptionally well-suited for large real-world datasets.

\subsection{Ablation Study}

In QIMF, query complexity measures efficiency, influenced by the number of shots $n_s$ for cost function evaluation, the number of Monte Carlo samples $n_b$, and the number of epochs $n_e$. The classical version of the mean field probabilistic model is characterized by a complexity of $O(n_w \cdot n_b \cdot n_e)$, whereas the dequantized version exhibits a reduced $O(n_s \cdot n_b \cdot n_e)$ complexity.

\begin{figure*}[ht]
    \centering
    \begin{minipage}{0.48\textwidth}
        \centering
        \includegraphics[width=\linewidth]{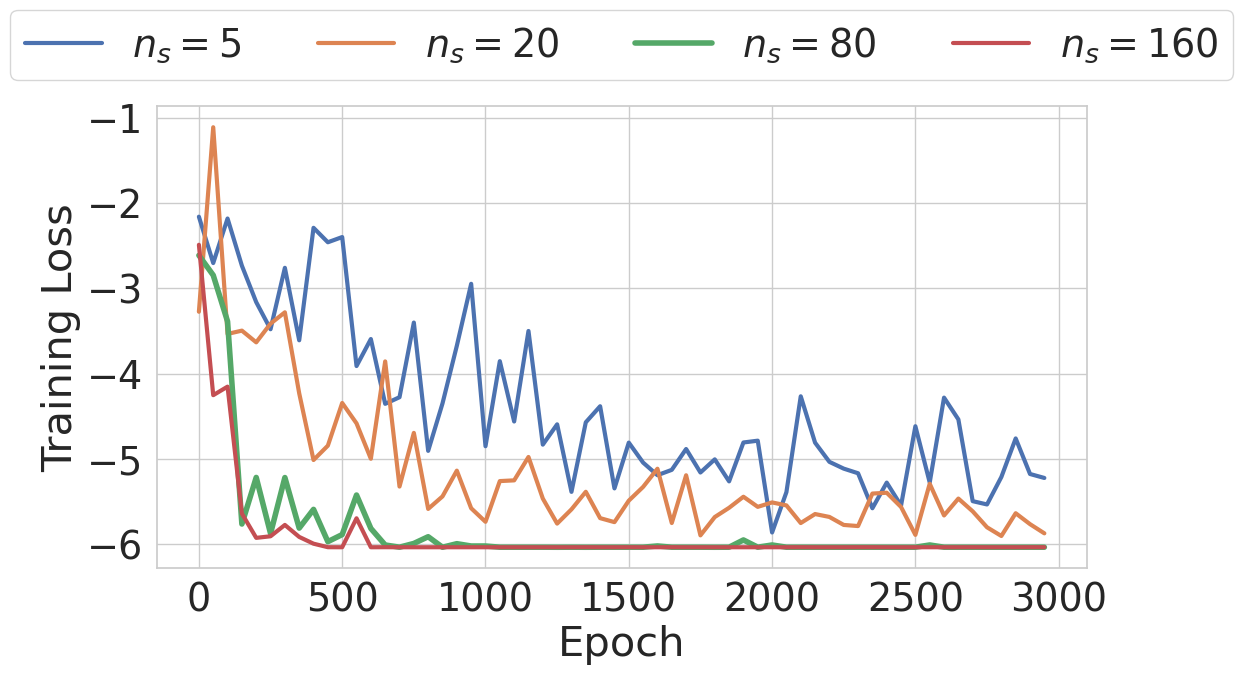}
        \caption{Comparison of different shot numbers $n_s$ on QIMF algorithm performance in solving QUBO problems using the WSBM with $N=5$ dataset.}
        \label{fig: case4}
    \end{minipage}\hfill
    \begin{minipage}{0.48\textwidth}
        \centering
        \includegraphics[width=\linewidth]{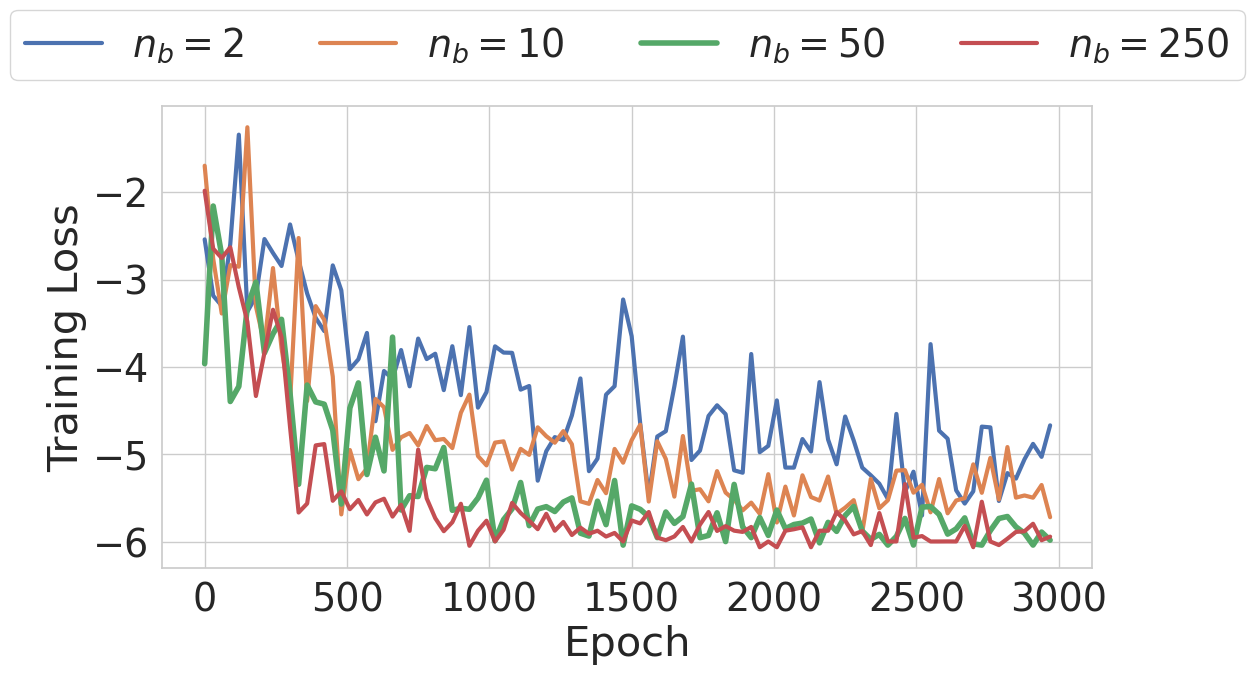}
        \caption{Comparison of different Monte Carlo samples $n_b$ on QIMF algorithm performance in solving QUBO problems using the WSBM with $N=5$ dataset.}
        \label{fig: case5}
    \end{minipage}
\end{figure*}

\textbf{Impact of Shot Count on QIMF Performance.}
We analyze the impact of different $n_s$ on the accuracy of QIMF and provide evidence on the significance of $n_s$ selection in the QUBO problems. We consider the WSBM with $N=5$ portfolio optimization problem, We set the number of Monte Carlo sampling for the QIMF algorithm at $n_b=100$. In Fig.~\ref{fig: case4}, we test a range of values from $n_s=5$ to $n_s=160$. For instance, with a relatively low $n_s$ such as $5, 20$, QIMF may not fully capture the underlying distribution of the problem space, leading to suboptimal solutions. Conversely, an excessively high $n_s=160$ results in redundant computations without significant gains in accuracy, thus diminishing the algorithm's speedup advantage. 
Thus, the optimal $n_s$ (80 in our case) must balance computational resource constraints and the accuracy demands of the specific problem, ensuring that the overall query complexity is optimized without compromising the solution quality.


\textbf{Impact of Monte Carlo Sample Count on QIMF Performance.}
Balancing the number of batches, $n_b$ is critical. We consider the WSBM with $N=5$ portfolio optimization problem, We set the number of shots for the QIMF algorithm at $n_s=20$. Fig.~\ref{fig: case5} shows that as $n_b$ increases from 2 to 250, the performance of the QIMF algorithm stabilizes, particularly for $n_b$ values from 50 to 250, where the performance levels are similar. However, smaller values of $n_b$ demonstrate greater efficiency, suggesting a trade-off between reducing instability and maximizing algorithmic efficiency. Thus, the optimal $n_b$ must strike a balance between minimizing instability for reliable outcomes and maintaining algorithmic efficiency.

\subsection{Scalability Analysis}\label{scal_a}


\begin{figure*}[ht]
    \centering
    \includegraphics[width=1\textwidth]{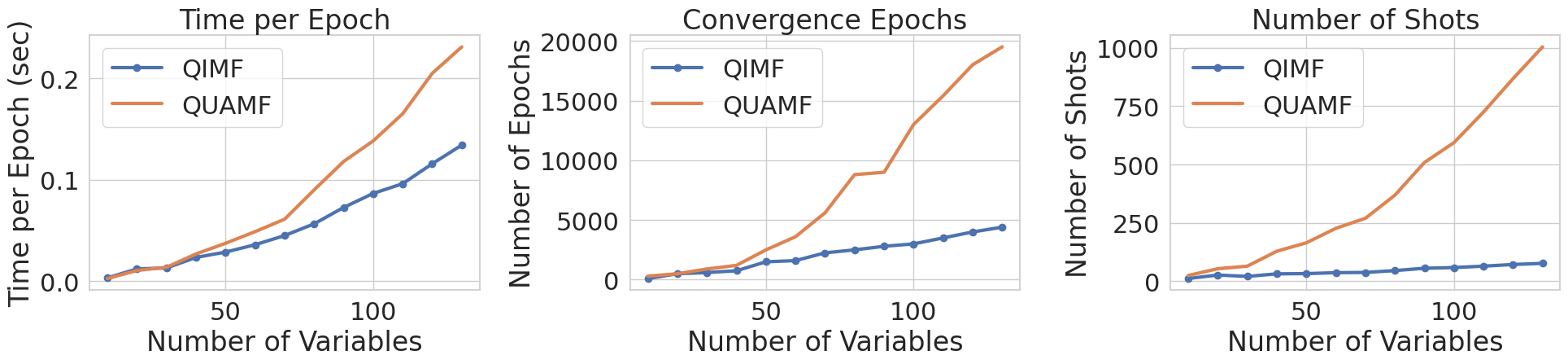}
    \caption{Scalability analysis of the QIMF and QUAMF, showing time per epoch as system size increases on the left, number of epochs required for convergence as system size increases in the middle, and number of shots needed as system size increases on the right.}
    \label{fig:scalability}
\end{figure*}


We employ the WSBM$(\mathbf{n}, N, P, W)$ datasets, where $N$ ranges from 10 to 130 and increases linearly with the number of variables. The matrix $P$ with diagonal elements at 0.2 and off-diagonal elements at 0.05. The matrix $W$ with diagonal entries following a normal distribution $N(0, 0.2)$ and off-diagonal entries as $N(0, 0.05)$.

To discuss the scalability of the QIMF model, it is essential to compare it with its classical counterpart. Here, QIMF stands for a mean field probabilistic model with a quantum-inspired cost function, whereas QUAMF represents a mean field probabilistic model with a classical cost function. This distinction is crucial as it highlights the impact of the cost function on the performance of both models.

The resource consumption of QIMF and QUAMF is analyzed through three key graphical representations, detailed in Figure~\ref{fig:scalability}. The first graph shows a slower increase in computational time for QIMF as the number of variables grows, compared to QUAMF. The second graph indicates a more gradual rise in the number of epochs needed for QIMF to converge, unlike the sharper increase for QUAMF. The third graph plots $n_s$ required by QIMF against that for QUAMF, demonstrating QIMF's more efficient resource usage. It should be noted that in the case of QUAMF, the corresponding $n_s$ is actually equal to $n_w$. These results indicate that the QIMF method offers significant advantages in resource consumption compared to QUAMF, particularly in larger-scale problems.

\section{Discussion}\label{discuss}

\begin{figure*}[t]
    \centering
    \begin{minipage}{0.48\textwidth}
        \includegraphics[width=\linewidth]{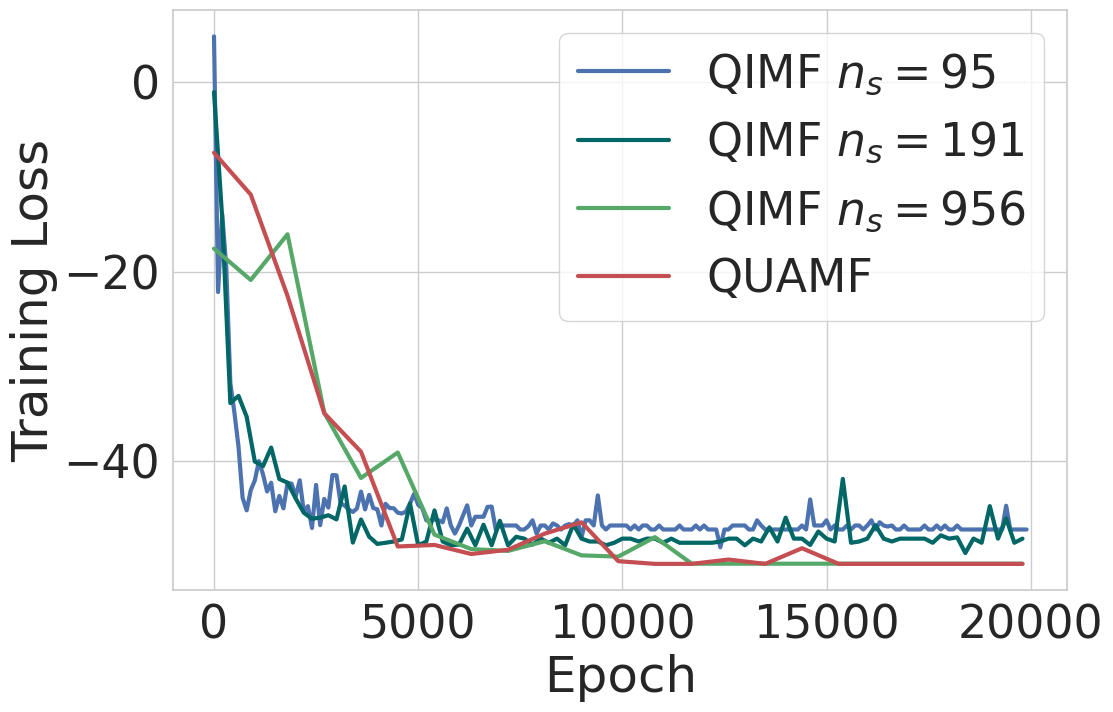}
        \caption{Performance comparison of QUAMF and QIMF under the WSBM dataset with $N=1$ for various \( n_s \) values.}
        \label{fig:Case6}
    \end{minipage}
    \begin{minipage}{0.48\textwidth}
        \includegraphics[width=\linewidth]{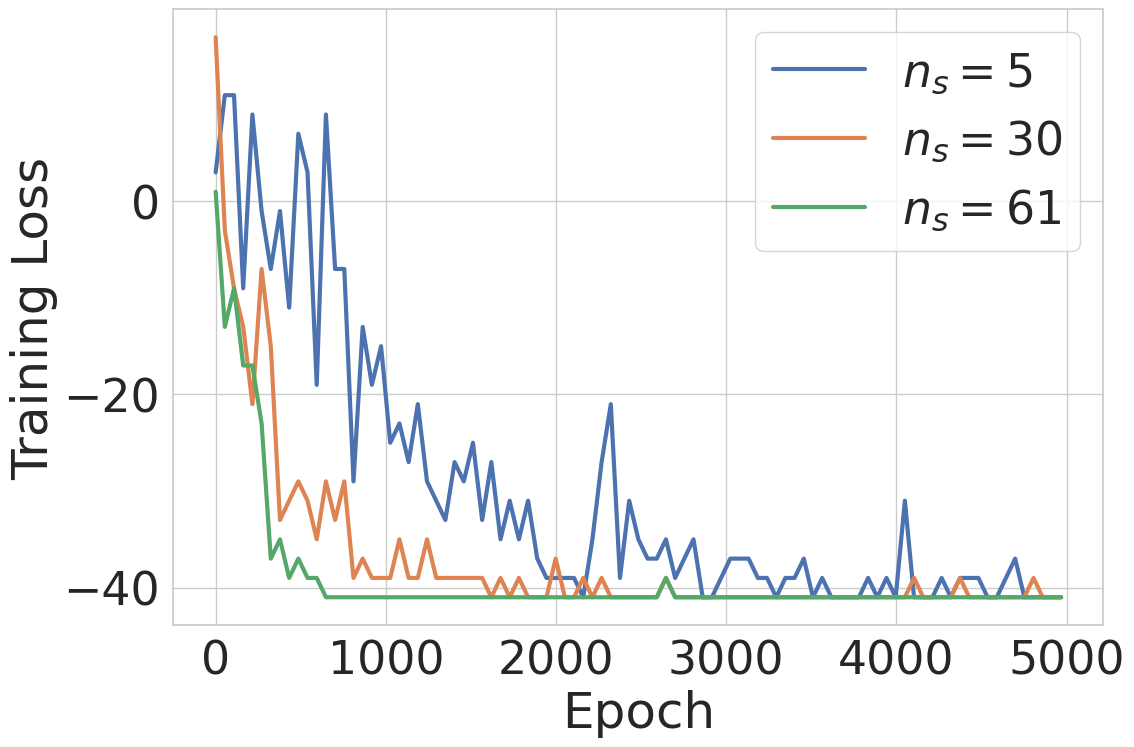}
        \caption{Performance comparison of QIMF under the unweighted SBM dataset with $N=5$ for various \( n_s \) values.}
        \label{fig:Case7}
    \end{minipage}   
\end{figure*}

\textbf{Theoretical Analysis.}
The number of shots $n_s$ is chosen to balance the computational efficiency and the accuracy in the QUBO problem. To select the number of shots in the QIMF method, we can leverage the structure of an approximate block diagonal matrix representing the $\text{WSBM}(\mathbf{n}, N, P, W)$. The assumptions are as follows:
\begin{itemize}
    \item Each element in vector $\mathbf{n}$ is equal and equivalent to $M$, and there are $N$ such elements in vector $\mathbf{n}$.
    \item The diagonal elements of matrix $P$ are $q$ and the off-diagonal elements are $p$, ensuring that the diagonal elements are significantly larger than the off-diagonal elements, which is typically represented as $q > p$.
    \item The diagonal parameters of matrix $W$ are $q'$ and the off-diagonal parameters are $p'$, ensuring that the diagonal elements are significantly larger than the off-diagonal elements, which is typically represented as $q' > p'$.
\end{itemize}

To capture the essential features of the WSBM, particularly the dense diagonal blocks, we aim to cover $n_d = M \times M \times N \times p$ elements. The total number of elements in the matrix is $n_w = M \times M \times N \times q + M \times M \times (N^2 - N) \times p$. The formula for the number of shots $n_s$ becomes $n_s = \frac{n_d}{n_w} \times n_w = \frac{1}{1 + (N-1) \frac{p}{q}} \times n_w$, simplifying to $n_s = \Omega\left(\frac{1}{N}\right) \times n_w$.

In the scenario where the dataset features an \(N\)-block structure with non-uniform distribution, our method demonstrates a computational speedup by a factor of \(N\). This speedup results from emphasizing the denser intra-block interactions, which have a more substantial influence on the system's state. 

When \(N=1\) or \(p = q\), the number of shots \(n_s\) closely approximates \(n_w\), indicating no speedup since the intra-block and inter-block interactions are equally probable and significant. In another case where \(p' = q'\), the weights assigned to intra- and inter-block interactions become indistinguishable. The probabilities $p$ and $q$ play a critical role in determining the accuracy of the results. In the following sections, we will discuss these two constrained scenarios in detail.

\textbf{Limitation.}
We evaluate performance under WSBM with $N=1$ dataset, with $n_w=956$ non-zero terms. We set $n_b=20$ and vary $n_s$ across 95, 191, and 956 values. The comparative performance of QUAMF and QIMF under these scenarios is illustrated in Fig.~\ref{fig:Case6}. At $n_s=95$ and $n_s=191$, the results show that the loss incurred by QIMF is higher than that by QUAMF because of the adverse effects of QIMF's shot reallocation strategy, which appears to overlook certain critical factors. As the number of shots increases to $n_s=956$, the loss value from QIMF aligns more closely with that of QUAMF. As there is convergence in performance at higher shot numbers, QIMF does not offer a computational advantage over QUAMF when operating under this WSBM setting with $N=1$.

In the amplitude-based allocation approach, If $a_m$ values in the cost function are nearly identical, the method struggles to achieve high-precision solutions when combined with the cost function scaled by $n_s/N$, indicating a lack of acceleration under these conditions. In our experiments on the Ising problem, we set $n_b=50$ and ensured that all terms represented in W are uniformly distributed. We then tested $n_s$ at three different values: 5, 30, and 61. The performance comparison of QIMF methods under various conditions is illustrated in Fig.~\ref{fig:Case7}. It is evident that for $n_s=61$, the method converges at about 800 epochs. When $n_s=30$, the convergence requires roughly twice as many epochs as it does for $n_s=61$. For $n_s=5$, the method does not converge even after 5,000 epochs, indicating that it does not offer a speedup in this case.

\section{Conclusion}
\label{conclusion}
We develop QIMF, a novel approach for addressing the QUBO problem through quantum-inspired techniques. Our method functions as a white box algorithm, utilizing grouping techniques derived from quantum measurement principles to dequantize cost function. We also use the amplitude-based shot allocation method, which refines the query complexity of the quantum-inspired cost function. The center optimizing strategy of our method is the mean field probabilistic model. Our theoretical analysis has confirmed that QIMF can indeed achieve polynomial speedup, not by relying on the main-solving strategy alone but through the effects of the allocated quantum-inspired cost function. 

Empirical validation of QIMF's effectiveness comes from four distinct case studies, each highlighting scenario capabilities. These studies demonstrate the method's substantial speedup, which we quantify as $\times N$, and its superior performance compared to heuristic algorithms, unsupervised learning models, and Monte Carlo sampling-based optimizers. In particular, our method outperforms the best heuristic algorithm by approximately $12.5\%$ in the S\&P500 2023 case. Additionally, we demonstrate the tradeoffs on choosing the number of shots $n_s$, and the number of samples $n_b$. Furthermore, the scalability of QIMF, compared to its classical counterpart QUAMF, demonstrates significant advantages, indirectly proving the efficacy of the quantum-inspired cost function. This ensures QIMF's applicability and efficiency in managing increasingly large datasets and complex optimization challenges, reinforcing its potential as a versatile tool in industrial applications.


\bibliographystyle{unsrt}
\bibliography{ref}

\begin{thebibliography}{10}

\bibitem{lang2022strategic}
Jonas Lang, Sebastian Zielinski, and Sebastian Feld.
\newblock Strategic portfolio optimization using simulated, digital, and quantum annealing.
\newblock {\em Applied Sciences}, 12(23):12288, 2022.

\bibitem{barahona1988application}
Francisco Barahona, Martin Gr{\"o}tschel, Michael J{\"u}nger, and Gerhard Reinelt.
\newblock An application of combinatorial optimization to statistical physics and circuit layout design.
\newblock {\em Operations Research}, 36(3):493--513, 1988.

\bibitem{lucas2014ising}
Andrew Lucas.
\newblock Ising formulations of many np problems.
\newblock {\em Frontiers in physics}, 2:74887, 2014.

\bibitem{kochenberger2014unconstrained}
Gary Kochenberger, Jin-Kao Hao, Fred Glover, Mark Lewis, Zhipeng L{\"u}, Haibo Wang, and Yang Wang.
\newblock The unconstrained binary quadratic programming problem: a survey.
\newblock {\em Journal of combinatorial optimization}, 28:58--81, 2014.

\bibitem{de1994cutting}
Caterina De~Simone and Giovanni Rinaldi.
\newblock A cutting plane algorithm for the max-cut problem.
\newblock {\em Optimization Methods and Software}, 3(1-3):195--214, 1994.

\bibitem{rendl2007branch}
Franz Rendl, Giovanni Rinaldi, and Angelika Wiegele.
\newblock A branch and bound algorithm for max-cut based on combining semidefinite and polyhedral relaxations.
\newblock In {\em Integer Programming and Combinatorial Optimization: 12th International IPCO Conference, Ithaca, NY, USA, June 25-27, 2007. Proceedings 12}, pages 295--309. Springer, 2007.

\bibitem{kagawa2020fully}
Hiroshi Kagawa, Yasuaki Ito, Koji Nakano, Ryota Yasudo, Yuya Kawamata, Ryota Katsuki, Yusuke Tabata, Takashi Yazane, and Kenichiro Hamano.
\newblock Fully-pipelined architecture for simulated annealing-based qubo solver on the fpga.
\newblock In {\em 2020 Eighth International Symposium on Computing and Networking (CANDAR)}, pages 39--48. IEEE, 2020.

\bibitem{kim2019comparison}
Yong-Hyuk Kim, Yourim Yoon, and Zong~Woo Geem.
\newblock A comparison study of harmony search and genetic algorithm for the max-cut problem.
\newblock {\em Swarm and evolutionary computation}, 44:130--135, 2019.

\bibitem{batsheva2024protes}
Anastasiia Batsheva, Andrei Chertkov, Gleb Ryzhakov, and Ivan Oseledets.
\newblock Protes: probabilistic optimization with tensor sampling.
\newblock {\em Advances in Neural Information Processing Systems}, 36, 2024.

\bibitem{nielsen2001quantum}
Michael~A Nielsen and Isaac~L Chuang.
\newblock {\em Quantum computation and quantum information}, volume~2.
\newblock Cambridge university press Cambridge, 2001.

\bibitem{pastorello2019quantum}
Davide Pastorello and Enrico Blanzieri.
\newblock Quantum annealing learning search for solving qubo problems.
\newblock {\em Quantum Information Processing}, 18(10):303, 2019.

\bibitem{farhi2014quantum}
Edward Farhi, Jeffrey Goldstone, and Sam Gutmann.
\newblock A quantum approximate optimization algorithm.
\newblock {\em arXiv preprint arXiv:1411.4028}, 2014.

\bibitem{zhang2022differentiable}
Shi-Xin Zhang, Chang-Yu Hsieh, Shengyu Zhang, and Hong Yao.
\newblock Differentiable quantum architecture search.
\newblock {\em Quantum Science and Technology}, 7(4):045023, 2022.

\bibitem{verteletskyi2020measurement}
Vladyslav Verteletskyi, Tzu-Ching Yen, and Artur~F Izmaylov.
\newblock Measurement optimization in the variational quantum eigensolver using a minimum clique cover.
\newblock {\em The Journal of chemical physics}, 152(12), 2020.

\bibitem{yen2020measuring}
Tzu-Ching Yen, Vladyslav Verteletskyi, and Artur~F Izmaylov.
\newblock Measuring all compatible operators in one series of single-qubit measurements using unitary transformations.
\newblock {\em Journal of chemical theory and computation}, 16(4):2400--2409, 2020.

\bibitem{zhu2023optimizing}
Linghua Zhu, Senwei Liang, Chao Yang, and Xiaosong Li.
\newblock Optimizing shot assignment in variational quantum eigensolver measurement.
\newblock {\em arXiv preprint arXiv:2307.06504}, 2023.

\bibitem{airoldi2008mixed}
Edo~M Airoldi, David Blei, Stephen Fienberg, and Eric Xing.
\newblock Mixed membership stochastic blockmodels.
\newblock {\em Advances in neural information processing systems}, 21, 2008.

\bibitem{decelle2011asymptotic}
Aurelien Decelle, Florent Krzakala, Cristopher Moore, and Lenka Zdeborov{\'a}.
\newblock Asymptotic analysis of the stochastic block model for modular networks and its algorithmic applications.
\newblock {\em Physical review E}, 84(6):066106, 2011.

\bibitem{karrer2011stochastic}
Brian Karrer and Mark~EJ Newman.
\newblock Stochastic blockmodels and community structure in networks.
\newblock {\em Physical review E}, 83(1):016107, 2011.

\bibitem{stanley2019stochastic}
Natalie Stanley, Thomas Bonacci, Roland Kwitt, Marc Niethammer, and Peter~J Mucha.
\newblock Stochastic block models with multiple continuous attributes.
\newblock {\em Applied Network Science}, 4:1--22, 2019.

\bibitem{faskowitz2018weighted}
Joshua Faskowitz, Xiaoran Yan, Xi-Nian Zuo, and Olaf Sporns.
\newblock Weighted stochastic block models of the human connectome across the life span.
\newblock {\em Scientific reports}, 8(1):12997, 2018.

\bibitem{xiao2019hybrid}
Yuchen Xiao and Ruzhe Zhong.
\newblock A hybrid recommendation algorithm based on weighted stochastic block model.
\newblock {\em arXiv preprint arXiv:1905.03192}, 2019.

\bibitem{yun2019optimal}
Se-Young Yun and Alexandre Prouti{\`e}re.
\newblock Optimal sampling and clustering in the stochastic block model.
\newblock {\em Advances in Neural Information Processing Systems}, 32, 2019.

\bibitem{peixoto2012entropy}
Tiago~P Peixoto.
\newblock Entropy of stochastic blockmodel ensembles.
\newblock {\em Physical Review E}, 85(5):056122, 2012.

\bibitem{caprara1999exact}
Alberto Caprara, David Pisinger, and Paolo Toth.
\newblock Exact solution of the quadratic knapsack problem.
\newblock {\em INFORMS Journal on Computing}, 11(2):125--137, 1999.

\bibitem{samorani2019clustering}
Michele Samorani, Yang Wang, Yang Wang, Zhipeng Lv, and Fred Glover.
\newblock Clustering-driven evolutionary algorithms: an application of path relinking to the quadratic unconstrained binary optimization problem.
\newblock {\em Journal of Heuristics}, 25:629--642, 2019.

\bibitem{izmaylov2019unitary}
Artur~F Izmaylov, Tzu-Ching Yen, Robert~A Lang, and Vladyslav Verteletskyi.
\newblock Unitary partitioning approach to the measurement problem in the variational quantum eigensolver method.
\newblock {\em Journal of chemical theory and computation}, 16(1):190--195, 2019.

\bibitem{Note1}
S\&P 500 stock data for the years 2022 and 2023 is sourced from~\protect \url {https://www.kaggle.com/datasets/pavankrishnanarne/s-and-p-500-stock-data-from-listing-day-to-2023/data}.

\bibitem{gurobi}
LLC Gurobi~Optimization.
\newblock Gurobi optimization, 2023.

\bibitem{zhang2023let}
Dinghuai Zhang, Hanjun Dai, Nikolay Malkin, Aaron Courville, Yoshua Bengio, and Ling Pan.
\newblock Let the flows tell: Solving graph combinatorial optimization problems with gflownets.
\newblock {\em Advances in neural information processing systems}, 2023.

\bibitem{vince2002framework}
Andrew Vince.
\newblock A framework for the greedy algorithm.
\newblock {\em Discrete Applied Mathematics}, 121(1-3):247--260, 2002.

\bibitem{bennet2021nevergrad}
Pauline Bennet, Carola Doerr, Antoine Moreau, Jeremy Rapin, Fabien Teytaud, and Olivier Teytaud.
\newblock Nevergrad: black-box optimization platform.
\newblock {\em ACM SIGEVOlution}, 14(1):8--15, 2021.

\bibitem{shor1994algorithms}
Peter~W Shor.
\newblock Algorithms for quantum computation: discrete logarithms and factoring.
\newblock In {\em Proceedings 35th annual symposium on foundations of computer science}, pages 124--134. Ieee, 1994.

\bibitem{grover1996fast}
Lov~K Grover.
\newblock A fast quantum mechanical algorithm for database search.
\newblock In {\em Proceedings of the twenty-eighth annual ACM symposium on Theory of computing}, pages 212--219, 1996.

\bibitem{kandala2017hardware}
Abhinav Kandala, Antonio Mezzacapo, Kristan Temme, Maika Takita, Markus Brink, Jerry~M Chow, and Jay~M Gambetta.
\newblock Hardware-efficient variational quantum eigensolver for small molecules and quantum magnets.
\newblock {\em nature}, 549(7671):242--246, 2017.

\bibitem{zhou2020quantum}
Leo Zhou, Sheng-Tao Wang, Soonwon Choi, Hannes Pichler, and Mikhail~D Lukin.
\newblock Quantum approximate optimization algorithm: Performance, mechanism, and implementation on near-term devices.
\newblock {\em Physical Review X}, 10(2):021067, 2020.

\end{thebibliography}

\appendix

\onecolumngrid

\section{Advanced Quantum Technique for WSBM-QUBO problem}

\subsection{Quantum Computing}\label{QC}
Quantum computing~\cite{nielsen2001quantum} is built on the principles of quantum mechanics, utilizing quantum bits (qubits), as its fundamental units. Unlike classical bits, which are strictly in states \(0\) or \(1\), qubits operate in states that can be described by superposition
$|\phi\rangle = a|0\rangle + b|1\rangle$, where \(a\) and \(b\) are complex numbers representing the probability amplitudes of the qubit being in states \(|0\rangle\) and \(|1\rangle\), respectively. To ensure the total probability of finding the qubit in either state sums to one, these amplitudes are subject to the normalization condition $|a|^2 + |b|^2 = 1$. This superposition allows qubits to process multiple possibilities at once, potentially increasing computation speed for some problems. Quantum entanglement is another key concept in quantum computing. It links pairs of qubits so that the state of one is connected to the state of others, regardless of distance. This unique correlation, absent in classical systems, doesn't allow for faster-than-light communication but enables quantum algorithms to outperform classical ones. Examples include Shor's algorithm~\cite{shor1994algorithms} for factoring integers and Grover's algorithm~\cite{grover1996fast} for searching databases more efficiently.

\subsection{Quantum Computing for WSBM-QUBO problem}\label{QA_QUBO}

\subsubsection{Weighted Stochastic Block Model}\label{WSBM}
The generalized WSBM significantly advances network analysis through its capability to model edge weights, accurately reflecting interaction strengths or capacities, which is essential for detailed analysis of complex network structures. The model is defined by a set number of nodes, $\mathbf{n}=\{n_1, n_2, \dots, n_N\}$, segmented into $N$ blocks, and utilizes a connectivity probability matrix $P \in \mathbb{R}^{N \times N}$. Each element \(P_{ij}\) within this matrix represents the probability of an edge existing between any node in block \(i\) and any node in block \(j\). Formally, it is defined as follows:
\begin{equation}
P = \begin{pmatrix}
P_{11} & P_{12} & \cdots & P_{1N} \\
P_{21} & P_{22} & \cdots & P_{2N} \\
\vdots & \vdots & \ddots & \vdots \\
P_{N1} & P_{N2} & \cdots & P_{NN}
\end{pmatrix},
\end{equation}
where \(P_{ij} \in [0, 1]\) for all \(i, j\). The accompanying weight distribution matrix $W$, which aligns with $P$, specifies the probabilistic distribution for the weights of these edges, thus enriching the model's detail and applicability. \(W_{ij}\) defines the probability distribution of the weights of edges that exist between nodes in block \(i\) and nodes in block \(j\). Depending on the application, various types of distributions can be used:
\begin{itemize}
\item If weights are normally distributed, \(W_{ij}\) could be represented as a normal distribution with parameters \(\mu_{ij}\) (mean) and \(\sigma_{ij}^2\) (variance), denoted as $\text{Dist}_{ij}(\theta_{ij})=\text{Norm}(\mu_{ij}, \sigma_{ij}^2)$.
\item If weights are exponential distributed, \(W_{ij}\) could be an exponential distribution with a rate parameter \(\lambda_{ij}\), denoted as $\text{Dist}_{ij}(\theta_{ij})=\text{Exp}(\lambda_{ij})$.
\item If weights are equal \(W_{ij}\) could be a uniform distribution between \(a_{ij}\) and \(b_{ij}\), expressed as $\text{Dist}_{ij}(\theta_{ij})=\text{Unif}(a_{ij}, b_{ij})$.
\end{itemize}
The overall weight distribution matrix $W$ is thus:
\begin{equation}
W = \begin{pmatrix}
\text{Dist}_{11}(\theta_{11}) & \text{Dist}_{12}(\theta_{12}) & \cdots & \text{Dist}_{1N}(\theta_{1N}) \\
\text{Dist}_{21}(\theta_{21}) & \text{Dist}_{22}(\theta_{22}) & \cdots & \text{Dist}_{2N}(\theta_{2N}) \\
\vdots & \vdots & \ddots & \vdots \\
\text{Dist}_{N1}(\theta_{N1}) & \text{Dist}_{N2}(\theta_{N2}) & \cdots & \text{Dist}_{NN}(\theta_{NN})
\end{pmatrix}, 
\end{equation}
where $\theta_{ij}$ specifies the statistical parameters governing how edge weights are distributed when an edge between respective blocks exists. In this work, we introduce a naming convention: WSBM-* problems. For example, when we apply this model to Quadratic Unconstrained Binary Optimization problems, we refer to them as WSBM-QUBO problems.

\subsubsection{Quantum Algorithms for WSBM-QUBO problem}
WSBM-QUBO problems can be efficiently encoded into the computational framework of quantum computers. This encoding involves the construction of a measurement operator \( H \) whose ground state corresponds to the optimal solution of the WSBM-QUBO problem. The measurement operator is generally expressed as a linear combination of tensor products of Pauli matrices $H = \sum_{i} a_i P_{1}^{(i)} \otimes P_{2}^{(i)} \otimes \cdots \otimes P_{n}^{(i)}$, where \( P_{j}^{(i)} \in \{I, X, Y, Z\} \) are Pauli operators acting on the \(j\)-th qubit and \( a_i \) are real coefficients.

The objective is to minimize the expectation value of this operator in a quantum state prepared by a parameterized quantum circuit $\min_{\theta} \langle \psi(\theta) | H | \psi(\theta) \rangle$, where \( |\psi(\theta)\rangle = U(\theta)|0\rangle \) is the state prepared from the all-zero state \( |0\rangle \) by the unitary operation \( U(\theta) \) parameterized by \( \theta \). Variational quantum eigensolver~\cite{kandala2017hardware} and quantum approximate optimization algorithm~\cite{zhou2020quantum} are notable quantum algorithms frequently employed for WSBM-QUBO problems.

\subsection{Efficient Measurement Strategies for Quantum Computing}\label{QWC}
To reduce the computational complexity involved in quantum measurements, the concept of qubit-wise commutativity is introduced. This method focuses on the commutativity of Pauli operators at each qubit position in multi-qubit systems, which simplifies the grouping of Hamiltonian terms for simultaneous measurement.

Consider two Pauli strings, \( H_i \) and \( H_j \), which are terms of the Hamiltonian \( H \). These terms are said to commute qubit-wise if and only if their corresponding Pauli operators at each qubit position commute. Specifically, \( H_i \) and \( H_j \) commute qubit-wise if \([H_{ik}, H_{jk}] = 0\) for every qubit \( k \), where \( H_{ik} \) and \( H_{jk} \) are the Pauli operators at the \( k \)-th position in \( H_i \) and \( H_j \) respectively.

This criterion is pivotal for the formation of QWC groups within the Hamiltonian. In such groups, every term \( H_i \) commutes qubit-wise with every other term in the same group. Conversely, terms in different groups do not satisfy the QWC condition. By organizing terms into these groups based on qubit-wise commutativity, it is possible to simplify the quantum measurement process, as terms within the same group can be measured simultaneously.

\section{Efficient Preprocessing of Quantum-Inspired Data}\label{Prep}

In quantum-inspired computing, effective data preprocessing is a key step in optimizing algorithm performance. Particularly before dequantization, precise preprocessing of the terms in the Hamiltonian can significantly impact the efficiency of the solution process and the quality of the results. In this section, we will detail how to efficiently preprocess the data items in the Hamiltonian.

Each term in the Hamiltonian can typically be described by a string containing quantum bit operators, such as 'Z' and 'I'. Specifically, we focus on those terms that contain only a single 'Z' operator. For these terms, we can determine their importance by comparing the absolute value of their coefficients to the sum of the absolute values of the coefficients of other strings that contain a 'Z' at the same position.

If the absolute value of the coefficient of a string containing only a single 'Z' exceeds the sum of the absolute values of the coefficients of all other strings with a 'Z' at the same position, then this coefficient can be considered diagonally dominant. Upon identifying diagonally dominant terms, we can directly determine the value of the corresponding position based on the sign of the coefficient: if positive, the solution of such variable is set to 0; if negative, it is set to 1.

This preprocessing method can effectively reduce the complexity of the problem, making the subsequent computational process more straightforward and efficient. However, it should be noted that although this method can reduce the computational burden, its resource overhead is \(O(n_w*\sum(\mathbf{n}))\). Therefore, depending on the scale of the specific problem and resource constraints, this preprocessing step is optional.

This efficient data preprocessing approach can reduce the complexity of problems, particularly in large-scale scenarios. However, as the complexity of large-scale data increases, the advantages of this preprocessing technique may not be significant. Therefore, the decision to apply this preprocessing method should be based on the specific requirements of the problem and the available resources.

\section{Derivation of Quantum Inspired Cost Function}
\label{QI_cost}

In this section, we derive the quantum-inspired cost function. Given the quantum cost function Eq.~\ref{quantum_cf}, the operators, denoted as $H_m$, take forms such as $Z_j$ or $Z_j Z_i$. Due to the commutative property, these operators satisfy
\begin{equation}
    [H_{i}, H_{j}] = 0,\label{commutator}
\end{equation}
indicating that they commute. Each operator $H_m$ can be represented in terms of its eigen decomposition
\begin{equation}
H_m = \sum_i \lambda_{mi} |\phi_i\rangle \langle\phi_i|,
\end{equation}
where $\lambda_{mi}$ are the eigenvalues, and $|\phi_i\rangle$ are the corresponding eigenvectors. We represent $H_m$ using tensor products of Pauli matrices
\begin{equation}
H_m = \bigotimes_{\otimes j} P_{mj}, \quad P_{mj} \in \{I, Z\},
\end{equation}
where $I$ is the identity matrix and $Z$ is the Pauli-Z matrix. Alternatively, the measurement operator can be expressed using projectors on the computational basis
\begin{equation}
H_m = \prod_{\otimes j} (|0\rangle \langle 0| + \eta_{mj} |1\rangle \langle 1|),
\end{equation}
where $\eta_{mj} \in \{-1, +1\}$. Here, $\eta_{mj} = +1$ corresponds to the identity matrix $I$, and $\eta_{mj} = -1$ corresponds to the Pauli-Z matrix $Z$. In the computational basis, written as
\begin{equation}
H_m = \sum_i \lambda_{mi} |\phi_i\rangle \langle \phi_i|,
\end{equation}
where $\phi_i \in \{0, 1\}^{\otimes n}$ represents the basis states across all qubits, and $\lambda_{mi} = \{-1, +1\}$. The cost function is then equal to
\begin{equation}
\text{Cost}(x) = \sum_m a_m \sum_i \lambda_{mi} |\langle \phi_i | x \rangle|^2,
\end{equation}
where $\phi_i$ and $x$ are classical variables representing quantum states. The cost function simplifies in cases where $\phi_i^\dagger x$ resolves to binary outcomes $\{0, 1\}$
\begin{equation}
\text{Cost}(x) = \sum_m a_m \lambda_{mi} I(\phi_i == x),
\end{equation}
where $I$ is the indicator function which equals 1 when $\phi_i$ equals $x$ and 0 otherwise.

\section{Derivation of Gradients of Model Parameters for Probabilistic Models}
\label{g_PM}

We begin with the expression for the objective function evaluated in Eq.~\ref{Eq: ell}. By applying the gradient $\nabla_{\alpha}$, we derive the following
\begin{align}
\nabla_{\alpha} \ell &= \sum_{x} \nabla_{\alpha} \left( \frac{P(x, \alpha)}{Z(\alpha)} \right) \text{Cost}(x) \\
&=\sum_x \frac{\nabla_{\alpha} P(x, \alpha)}{Z(\alpha)} \text{Cost}(x) - \frac{\nabla_{\alpha} Z(\alpha)}{Z(\alpha)}  P(x, \alpha) \text{Cost}(x)\\
&= \sum_{x} \frac{P(x, \alpha)}{Z(\alpha)} \left( \frac{\nabla_{\alpha} P(x, \alpha)}{P(x, \alpha)} \text{Cost}(x) \right) - \frac{P(x, \alpha)}{Z(\alpha)} \left( \frac{\nabla_{\alpha} Z(\alpha)}{Z(\alpha)} \text{Cost}(x) \right) \\
&= \sum_{x} \frac{P(x, \alpha)}{Z(\alpha)} \left(  \frac{\nabla_{\alpha} P(x, \alpha)}{P(x, \alpha)} \text{Cost}(x) \right) - \left(\sum_{x} \frac{P(x, \alpha)}{Z(\alpha)} \text{Cost}(x)\right) \left(\sum_{x} \frac{P(x, \alpha)}{Z(\alpha)} \nabla_{\alpha} P(x, \alpha) \right) \\
&= \sum_{x \in P} \nabla_{\alpha} \ln P(x, \alpha) \text{Cost}(x) - \left(\sum_{x \in P} \nabla_{\alpha} \ln P(x, \alpha)\right) \left(\sum_{x \in P} \text{Cost}(x)\right).
\end{align}

Given the property of normalized probability distributions, we have $\sum_{x \in P} \nabla_{\alpha} \ln P(x, \alpha) = 0$. Therefore, the gradient of the objective function simplifies to
\begin{equation}
\nabla_{\alpha} \ell = \sum_{x \in P} \nabla_{\alpha} \ln P(x, \alpha) \text{Cost}(x).
\end{equation}

To further compute the gradient $\nabla_{\alpha_{ij}} \ln P(x_i = m, \alpha)$, consider the following derivation
\begin{align}
\nabla_{\alpha_{ij}} \ln P(x_i = m, \alpha) &= \nabla_{\alpha_{ij}} \ln \frac{e^{\alpha_{im}}}{\sum_k e^{\alpha_{ik}}} \\
&= \nabla_{\alpha_{ij}} \left(\alpha_{im} - \ln \left(\sum_k e^{\alpha_{ik}}\right)\right) \\
&= \delta(j = m) - \frac{e^{\alpha_{ij}}}{\sum_k e^{\alpha_{ik}}} \\
&= \delta(j = m) - P(x_i = m),
\end{align}
where $\delta(j = m)$ is the Kronecker delta, which equals 1 when $j = m$ and 0 otherwise.

\section{Optimal Block Size Selection N for Portfolio Optimization}\label{PO}

Industry sectors play a crucial role in the classification of S\&P 500 stock datasets. In this case, we analyzed a subset of 30 portfolios spanning 22 distinct industries, utilizing data from June 2023.

\begin{figure*}[htbp]
    \centering

    \begin{minipage}{0.4\textwidth}
        \includegraphics[width=\linewidth]{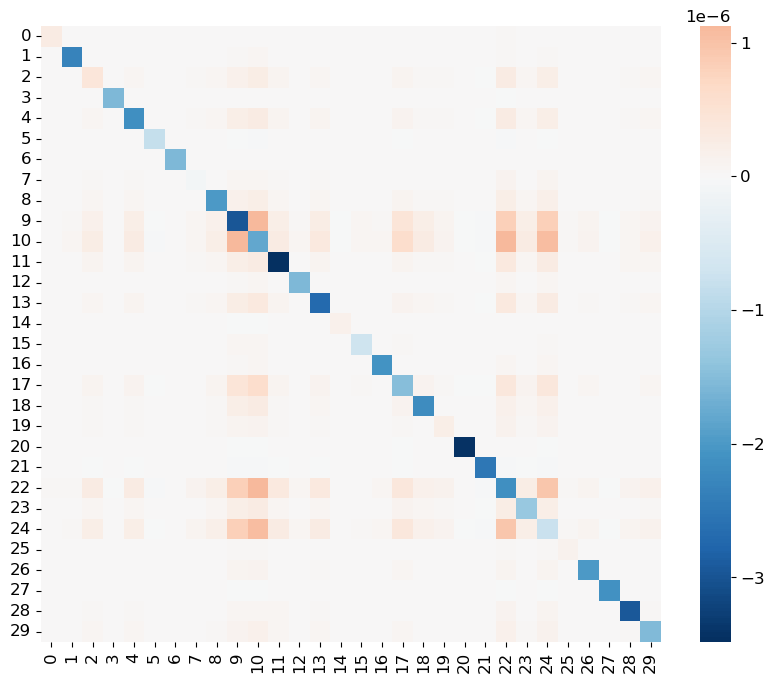}
        \caption{Heatmap displaying the problem distribution across 22 industries in the subset of 30 portfolios.}
        \label{fig:SP30_1}
    \end{minipage}
    \hfill 
    \begin{minipage}{0.58\textwidth}
        \includegraphics[width=\linewidth]{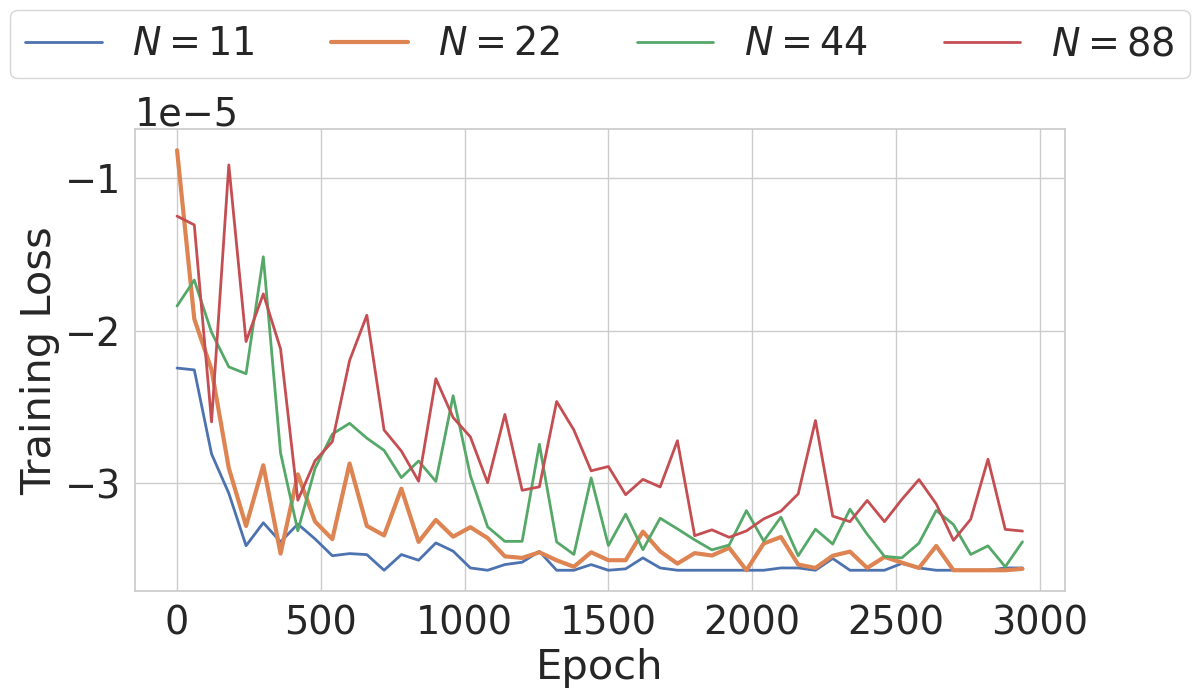}
        \caption{Comparison of training losses within the different N.}
        \label{fig:SP30_2}
    \end{minipage}
    
\end{figure*}

The heatmap depicted in Figure~\ref{fig:SP30_1} illustrates the distribution of data related to risk minus return in our portfolio optimization model. Specifically, it visualizes the calculated values of $\lambda x^\dagger V x- (1-\lambda)x^\dagger r$, where $\lambda$ represents the trade-off between risk (modeled by the variance-covariance matrix $V$) and return ($r$). As $x$ is constrained to be a binary vector in $\{0, 1\}^n$, the objective function transforms into $x^\dagger(\lambda V - (1-\lambda) r) x$. Thus, the distribution of data is $\lambda V - (1-\lambda) r$.

Additionally, Fig.~\ref{fig:SP30_2} compares the training losses for different block sizes $N=11$, $N=22$, $N=44$, and $N=88$. The results indicate that the performance for $N=11$ is similar to that for $N=22$. However, the performance deteriorates significantly for $N=44$ and $N=88$. Considering both efficiency and performance, $N=22$ appears to be the most suitable choice, though it is not the optimal one, as the heatmap does not demonstrate optimal blocking. To achieve optimal classification, a more fine-grained analysis would be required, which would entail additional resources. In this paper, we utilize a coarse-grained classification approach based on the number of industries.

\end{document}